\documentclass[preprint]{elsarticle}

\usepackage{float}

\usepackage{lineno,hyperref}
\modulolinenumbers[5]

\usepackage{amsmath}
\usepackage{amsfonts}
\usepackage{amssymb}
\usepackage{bm}
\allowdisplaybreaks[4]

\usepackage{tikz}

\usepackage{url}

\usepackage{multirow}

\usepackage{threeparttable}

\usepackage{mdwlist}

\usepackage{pgfplots}
\pgfplotsset{width=6cm,compat=1.13}

\usepackage{color}

\newtheorem{definition}{Definition}

\journal{Renewable and Sustainable Energy Reviews}

\begin{document}

\begin{frontmatter}

\title{Grid-side Flexibility of Power Systems  in Integrating Large-scale Renewable Generations: A Critical Review on Concepts, Formulations and Solution Approaches}

\author[mymainaddress]{Jia Li}

\author[mymainaddress]{Feng Liu\corref{mycorrespondingauthor}}
\cortext[mycorrespondingauthor]{Corresponding author}
\ead{lfeng@tsinghua.edu.cn}

\author[mysecondaryaddress]{Zuyi Li}

\author[mytertiaryaddress]{Chengcheng Shao}

\author[myfourthaddress]{Xinyuan Liu}

\address[mymainaddress]{State Key Laboratory of Power System, Department of Electrical Engineering, Tsinghua University, Beijing 100084, China}
\address[mysecondaryaddress]{Robert W. Galvin Center for Electricity Innovation at Illinois Institute of Technology, Chicago, IL 60616 USA}
\address[mytertiaryaddress]{School of Electrical Engineering and State Key Laboratory on Electrical Insulation and Power Equipment in Xi’an Jiaotong University, Xi'an, 710049, China.}
\address[myfourthaddress]{Shanxi Electric Power Research Institute, China.}

\begin{abstract}
Though considerable effort has been devoted to exploiting generation-side and demand-side operational flexibility in order to cope with uncertain renewable  generations, grid-side operational flexibility has not been fully investigated. In this review, we define grid-side flexibility as the ability of a power network to deploy its flexibility resources to cope with the changes of power system state, particularly due to variation of renewable  generation. Starting with a survey on the metrics of operational flexibility, we explain the definition from both physical and mathematical point of views. Then conceptual examples are presented to demonstrate the impacts of grid-side flexibility graphically, providing a geometric interpretation for a better understanding of the concepts. Afterwards the formulations and solution approaches in terms of grid-side flexibility in power system operation and planning are reviewed, based on which  future research directions and challenges are outlined.
\end{abstract}

\begin{keyword}
{Operational flexibility}\sep renewable generation integration\sep uncertainty\sep
{robust optimization}\sep
{stochastic optimization}
\end{keyword}

\end{frontmatter}


\section{Introduction}\label{sec_introduction}

With the rapid growing penetration of renewable generation, significant challenges of security and economics have arisen in power system operation due to the intermittent and stochastic nature and low predicability of renewable generation. Thus sufficient power system flexibility is required to cope with the new issues which have not been experienced before.
This paper is focused on the operational flexibility on the grid side of power systems with regard to operation and planning. For short, we call it ``grid-side flexibility''. The aim of this paper is to provide a comprehensive understanding on its concepts, formulations and solution approaches. Based on the review, the future research directions and challenges in optimally utilizing grid-side flexibility to facilitate a secure operation of the power system with high-penetration renewable generation can be figured out.

Power system flexibility is not a new concept, as power systems have always had to utilize generation resources, control systems and business practices to ensure that system supply-demand balance can be retained within the industry standards \cite{Kehler2011a}. Conventional methods to accommodate load uncertainty include regulating reserve, automatic generation control (AGC) and so on. However, these methods may not be able to provide sufficient flexibility to address the inherent uncertainty and volatility of renewable energy generation, which cannot be forecasted as accurately as electricity demand nowadays. To cope with the great challenge, new technologies have been proposed and constantly developed.

The flexibility of power system can be generally divided into three categories: generation side, grid side, and demand side. In the generation side, different kinds of approaches have been applied in unit commitment and economic dispatch to enhance generation-side flexibility. Stochastic optimization has been studied extensively, which explicitly incorporates uncertainty in the decision process \cite{Wang2008,Pappala2009,Zhang2011}. Most of the models rely on pre-sampling discrete scenarios and aim to minimize the expected cost. Instead of scenarios, interval optimization uses confidence intervals to characterize uncertainty, and derives optimistic and pessimistic solutions for satisfying system operational requirements \cite{Wu2012,Wang2011}. Taking advantage of uncertainty set, a probabilistic distribution is not required in robust optimization \cite{Jiang2012a,Bertsimas2013,Ye2015c}.
 An optimal solution is obtained, which immunizes against all the uncertain data contained within the given uncertainty set. Other approaches, such as fuzzy mixed integer program \cite{Daneshi2009a} and minimax regret program \cite{Jiang2013a}, have also been applied to help accommodate large-scale volatile renewable  generation.

In the demand side, demand response management is acknowledged to improve the operational flexibility \cite{Aghaei2013,Oconnell2014}. It has been applied in different aspects, such as unit commitment \cite{Zhao2013,DeJonghe2014}, real-time dispatch \cite{Madaeni2013}, and regulation provision \cite{Zhang2014} to cope with uncertainties in power system operation. Different schemes have been suggested, such as price-based \cite{Li2014} and coupon-based mechanisms \cite{Fang2015} to name a few, in order to create  incentives and encourage participation of demand response.

Compared with the flexibility from  generation and  demand sides, the grid-side flexibility has drawn less attention to date.  Physically, power network  provides sufficient capacity for transferring power from generation plants to consumers, which is traditionally considered as a fixed structure.
 However, due to the integration of large-scale renewable generations, the operation of power network has been pushed much closer to its technical limits than before. As a consequence, new problems arise which can hardly be resolved when there is lack of grid-side flexibility. First, the actual outputs of uncertain renewable generation may remarkably deviate from the forecast, causing  high operational risk or  congestion. The congestion can further increase operating cost and limit the usage of available flexibility resource from the generation and/or demand sides. Second, renewable generation usually requires the support of sufficient reactive power, imposing additional risk of voltage instability on power system operation.

In this context, novel technologies have been developed to exploit the potential flexibility of power network, making it a crucial supplement to generation-side and demand-side flexibility, and also an effective approach to address issues associated with congestion and voltage stability. Line switching has changed the traditional idea of a power network with a fixed topology, enabling power flow control by switching lines. Flexible AC transmission systems (FACTS) and high-voltage direct current (HVDC) technologies, have introduced more controllability into transmission networks \cite{Singh2015,Gandoman2018,Korompili2016,Bianchi2016}. They are mainly utilized to underpin reactive power compensation, voltage control, and power flow control \cite{Zhang2012a}, for the sake of increasing transmission capacity and  power system security. Due to the powerful controllable power electronic devices, much faster controllability can be provided in comparison to the  generation-side and demand-side resources.

In the literature, some reviews have focused on power system flexibility. Reference \cite{Kondziella2016} classifies the scientific approaches, that have been used in flexibility demand studies, into technical, economic, and market potential categories, based on the results from German and European energy systems. Reference \cite{Lund2015} reviews both supply-side and demand-side approaches, technologies, and strategies to enable high levels of variable renewable energy. Reference \cite{Alizadeh2016} classifies and discusses the possible flexibility impacts, including super short-term, short-term, mid-term, and long-term. It has been pointed out that new transmission technologies can enhance grid-side flexibility \cite{Alizadeh2016}. To the best of our knowledge, there is still lack of comprehensive reviews on grid-side flexibility so far.

\begin{figure}[htb]
\centering
\includegraphics[width=4in]{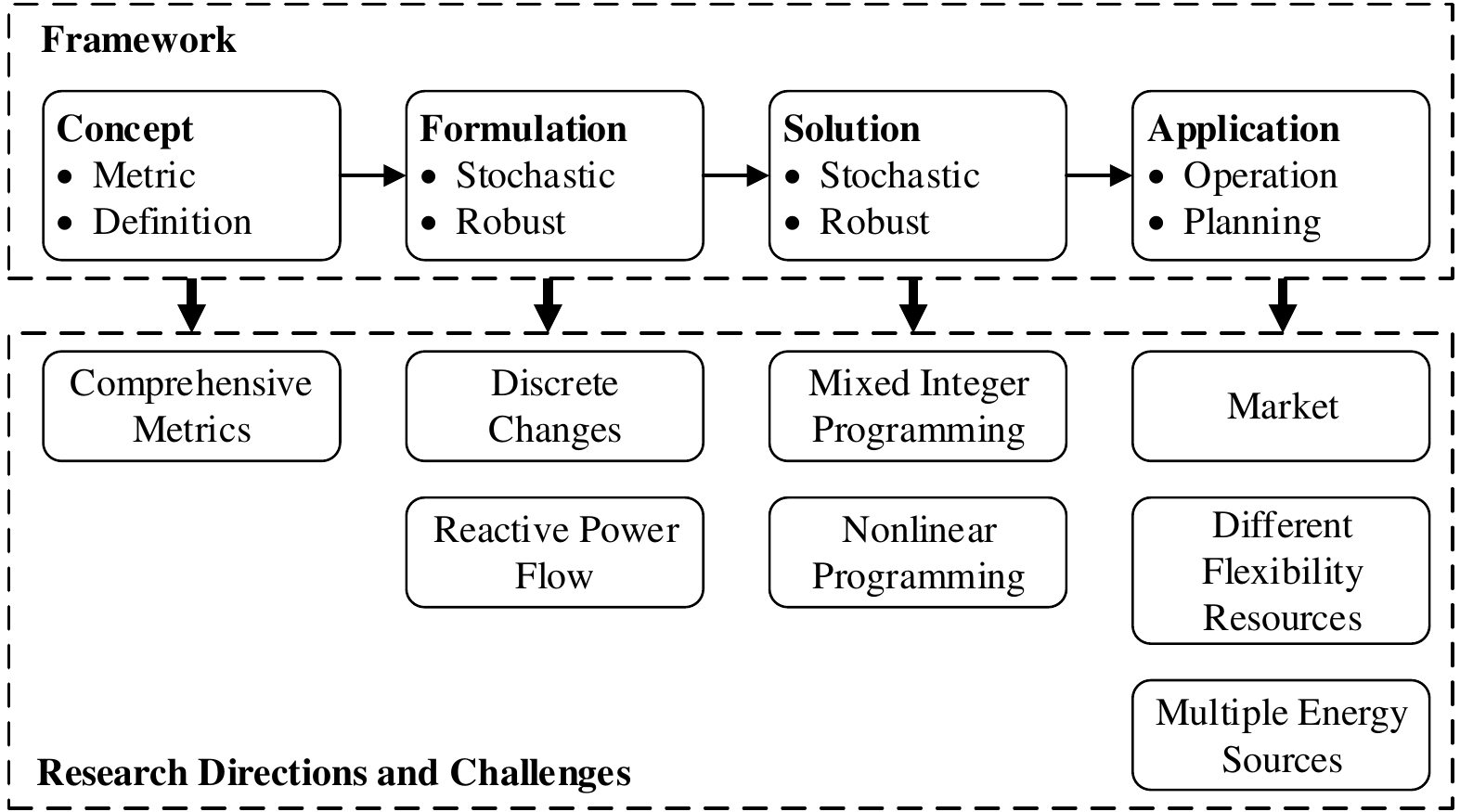}
\caption{Structure of the review.}
\label{fig_framework}
\end{figure}

This paper intends to provide a comprehensive understanding on ``grid-side flexibility'' of power systems, including its concepts, formulations as well as solution approaches, particularly when large-scale volatile renewable generation is integrated.  
To this end, the concept of grid-side flexibility is introduced first. Then the physical meaning of grid-side flexibility is explained, with the mathematical formulation of grid-side flexibility region. Moreover, geometric interpretation is explained via straightforward visualization, providing an intuitive understanding on the effects of various resources of grid-side flexibility. It also reveals the potential benefits of grid-side flexibility in accommodating uncertainty in power systems. For making better use of grid-side flexibility, the state-of-the-art studies on the theory and application are reviewed, including problem formulations and solution approaches. Last but not least, future research directions and challenging problems are outlined. Figure \ref{fig_framework} summarizes the structure of the paper.

The rest of the paper is organized as follows. Section II introduces the concept of grid-side flexibility. Section III reviews the formulations related to grid-side flexibility in power system operation and planning. Section IV reviews the associated solution approaches. Section V discusses the future research directions and challenging problems. Final remarks about the literature review and outlook of grid-side flexibility close the paper. 

\section{Concepts of Grid-side Flexibility}
In this section, first we review the metrics of power system flexibility. Then we explain our definition of grid-side flexibility in both physical and mathematical point of views. Moreover, geometric interpretation is explained via straightforward visualization, providing an intuitive understanding on the impacts of different grid-side flexibility resources.

\subsection{Literature Review of Flexibility Metrics}
In the literature, a number of metrics have been proposed to quantify power system flexibility with regard to operation or planning. However, a comprehensive metric which explicitly quantifies grid-side flexibility with regard to renewable energy generation has not been proposed so far. As a supplement to generation-side and demand-side flexibility, grid-side flexibility is determined not only by the topology and parameters of power networks, but also associated with the constraints of generation-side and demand-side resources.

\subsubsection{Power System Operation}
In power system operation, the metrics of flexibility mainly focus on generation-side flexibility, without consideration of grid-side flexibility. In other words, power networks are considered with fixed topologies and parameters.

Some metrics are based on the parameters of generators. In \cite{Menemenlis2011}, a flexibility index borrowed from the process control literature is proposed to evaluate an operation strategy which provides balancing reserves to mitigate wind power generation uncertainty. In \cite{Ma2013}, a metric is presented to quantify the ability of a generator to cope with the flexibility requirement as below.
\begin{equation}\label{eq_metric1}
	flex(i) = \frac{\frac{1}{2}\left[P_{\max}(i) - P_{\min}(i)\right] + \frac{1}{2}\left[Ramp(i)\Delta t\right]}{P_{\max}(i)}
\end{equation}
where $P_{\max}(i)$ and $P_{\min}(i)$ are the maximum capacity and the minimum stable generation of conventional generator $i$, while $Ramp(i)$ is the average value of the ramping up and down rates of generator $i$ per time period. $\Delta t$  stands for one time period, say, one hour. This index is further extended in \cite{VanHorn2014} for evaluating generation-side and demand-side flexibility while taking into account the impacts of transmission network. Inspired by the commonly used reliability metric, i.e., loss-of-load probability (LOLP), the lack-of-ramp probability (LORP) is proposed for  real-time economic dispatch \cite{Thatte2015a}. The system-wide LORP for the ramp-up case is defined in (\ref{eq_metric2}), which provides an assessment of the adequacy of the available system ramping capability from dispatched generators to meet both expected changes and uncertainty in forecasted net load \cite{Thatte2014}.
\begin{equation}\label{eq_metric2}
	LORP_s^{up,\tau}[t] = Pr\left(\sum_{i\in I}\left\{P_i^g[t] + \min\left(\tau R_i,P_i^{\max} - P_i^g[t]\right)\right\} < \widetilde{P}_s^l[t+\tau]\right)	
\end{equation}	 
where $P_i^g[t]$ is the dispatched output of generator $i$ at time $t$. $R_i$  is the one-interval ramp rate, and $P_i^{\max}$ is the maximum output of generator $i$. $\widetilde{P}_s^l[t+\tau]$ denotes the system-wide net load for the interval $\tau$ time steps in the future.

Stemming from conventional wisdom which only adopts one dimension of parameters, a few approaches use multiple dimensions of parameters to construct a multi-dimensional image for flexibility quantification. In \cite{Lu2010}, the concept of ``flying-brick'' is proposed, including three aspects: capacity, ramp rate, and energy. The metric is visualized as a three-dimensional probability box at each scheduling hour. Another three-dimensional space where the coordinates are magnitude, ramp rate, and ramp duration is presented in \cite{Dvorkin2014a} to characterize the intra-hour net load deviations where flexibility deployment is required. Based on the metrics of power provision capacity, power ramp-rate capacity, and energy provision capacity, a flexibility cube \cite{Ulbig2014} is constructed for quantifying and visualizing the technically available operational flexibility. The flexibility envelop concept follows the same line but considers the time-evolution of reserve requirements, which resembles the shape of an envelope or a cone \cite{Nosair2015}. The flexibility envelope concept is extended by characterizing operating reserve dynamically in \cite{Nosair2015a}.

Built on the existence of feasible operation strategies, operational flexibility is quantified in alternative way. In \cite{Bucher2015a} and \cite{Bucher2015}, operational flexibility is defined as the ability of a system to react to a disturbance sufficiently fast in order to keep the system secure. Focused on the concept of locational flexibility,  \cite{Bucher2015a} defines the ability of a system to contain a certain disturbance at a given node in terms of ramping rate, power and energy. Defined on the space spanned by the tie-line power injections, the flexibility sets \cite{Bucher2015} are presented to characterize the available operational flexibility in a multi-area power system, where N-1 security criterion is also taken into account. The deviation of external sources $p_e$ together with a corresponding reaction of internal sources $p_i$ that are feasible form the active flexibility set as in \eqref{eq_metric3}.
\begin{equation}\label{eq_metric3}
	F = \left\{(p_i,p_e)\in \mathbb{R}^{n_i}\times\mathbb{R}^{n_e} | 
	\begin{bmatrix}
		C_{i,N} \\
		C_{i,N-1}
	\end{bmatrix} p_i + 
	\begin{bmatrix}
		C_{i,N} \\
		C_{i,N-1}
	\end{bmatrix} p_e \le
	\begin{bmatrix}
		b_N \\
		b_{N-1}
	\end{bmatrix}
	\right\}	  	
\end{equation}
where the operational constraints are written in a matrix form. Subscript $``N"$ stands for a normal operation, while $``N-1"$ denotes N-1 security criterion. Inspired by the spirit of robust optimization, several operational flexibility metrics are constructed indicating the largest feasible operating regions of renewable generation without sacrificing power system security \cite{Zhao2015,Zhao2016,WangC.LiuF.WangJ.WeiW.Mei2016}. Similarly, dispatchability is defined as the largest set in the wind power generation uncertainty space such that the second stage dispatch is feasible if the uncertain data does not exceed its boundaries \cite{Wei2015a,Wei2015e}.

\subsubsection{Power System Planning}

In generation expansion planning, the insufficient ramping resource expectation (IRRE) metric is proposed to measure operational flexibility, which is derived from traditional generation adequacy metrics and is focused on generation-side flexibility \cite{Lannoye2012e}. Calculation of the IRRE follows a similar structure to the renowned reliability metric, loss-of-load expectation (LOLE). A distribution of the available flexible resources is formed for each direction in individual time horizon. Then the probability that the system experiences insufficient ramp resources at each observation $t$, over each time horizon $i$ and direction, is calculated from the available flexibility distribution. Thus, the overall metric is computed as shown in \eqref{eq_metric4}.
\begin{equation}\label{eq_metric4}
	IRRE_{i,+/-} = \sum_{\forall t\in T_{+/-}}AFD_{i,+/-}(NLR_{t,i,+/-} - 1)
\end{equation}
where $NLR_{t,i,+/-}$ denotes the net load ramp at observation $t$ in either
direction, and $AFD_{i,+/-}(X)$ indicates the probability that $X$ MW, or less, of flexible resource will be available. Since the temporal correlation between the flexibility available and the flexibility required is broken in this process, the unique value of IRRE is to highlight the time horizons where the system may have insufficient flexibility to meet unexpected changes of net load, instead of representing the precise number of expected periods of insufficient flexibility \cite{Lannoye2012d}. A high-level methodology to determine IRRE without requiring considerable off-line data and computation is presented in \cite{Lannoye2012} for planning. The index of periods of flexibility deficit (PFD) is also introduced to quantify the balance between the flexibility available and the flexibility required in a deterministic manner. The IRRE and PFD metrics can be used to assess the operational flexibility of a power system considering transmission network constraints \cite{Lannoye2015}. The results indicate that transmission network has significant effects on the usage of operational flexibility.

In transmission expansion planning, a few metrics have been proposed considering uncertainties in load growth \cite{Lu2005} and generation expansion \cite{Bresesti2003}. In \cite{Lu2005}, transmission expansion flexibility is quantified as a function of expected unserved load due to branches contingency as below.
\begin{align}
	TEF &= \int_{medium\_load}^{high\_load}g(l)dl \\
	g(l) &= \sum_k\sum_jL_{kj}P_j 
\end{align}
where $l$ is the customer required demand. $high\_load$ and $medium\_load$ are the values of peak load in high-growth and medium-growth scenarios, respectively. $L_{kj}$ represents the load shedding at bus $k$ to alleviate line overloads arising due to contingency $j$. $P_j$ denotes the probability of the occurrence of outage $j$. $g(l)$ is the expected unserved power. A more flexible plan is defined as the one that offers better adjustability for updating when load level increases. In \cite{Bresesti2003}, the flexibility of transmission system is defined as the attitude of the transmission system to retain a desired standard of reliability, at reasonable operation costs, when the generation scenarios change. Network flexibility is evaluated with regard to uncertainties of generation expansion and operation in a free electricity market. The proposed flexibility metric is an average of distribution factors of branches, weighted with the current margins of the corresponding branches as shown in \eqref{eq_metric7}. 
\begin{equation}\label{eq_metric7}
	\overline{C}_{ij}(I_{marg}) = \frac{\sum_{i=1,j=1}^{N_b}I_{marg}^{ij}C_{tot}^{ij}}{\sum_{i=1,j=1}^{N_b}I_{marg}^{ij}}(i\ne j)
\end{equation}
where, $N_b$ is the number of nodes; $I_{marg}$ is the current margin on branch $ij$; $C_{tot}^{ij}$ is the total distribution factor. A power network with the smaller value of the metric is considered more flexible, having more capability to deal with changes in power flow, and being less affected by variable generation. Two indices, namely technical uncertainty scenarios flexibility index and technical economical uncertainty scenarios flexibility index, are extended from \cite{Bresesti2003}. One of them consists of power flow margins alone, and the other also includes economic terms \cite{Capasso2003}.

Table \ref{tab_metric} summaries the flexibility metrics in the literature.

\begin{table}[htbp]
  \centering
  \footnotesize
  \caption{Flexibility Metrics in the Literature}
    \begin{tabular}{ccccc}
    \hline
    Ref. No. & \multicolumn{2}{c}{Research Aspect} & \multicolumn{2}{c}{Flexibility Resource} \\
          & Operation & Planning & Generation & Transmission \\
    \hline
	\cite{Menemenlis2011}               & $\surd$     &      & $\surd$     &  \\
	\cite{Ma2013}                       & $\surd$     &      & $\surd$     &  \\
	\cite{VanHorn2014}                  & $\surd$     &      & $\surd$     &  \\
	\cite{Thatte2015a}                  & $\surd$     &      & $\surd$     &  \\
	\cite{Thatte2014}                   & $\surd$     &      & $\surd$     &  \\
	\cite{Lu2010}                       & $\surd$     &      & $\surd$     &  \\
	\cite{Dvorkin2014a}                 & $\surd$     &      & $\surd$     &  \\
	\cite{Ulbig2014}                    & $\surd$     &      & $\surd$     &  \\
	\cite{Nosair2015}                   & $\surd$     &      & $\surd$     &  \\
	\cite{Nosair2015a}                  & $\surd$     &      & $\surd$     &  \\
	\cite{Bucher2015a}                  & $\surd$     &      & $\surd$     &  \\
	\cite{Bucher2015}                   & $\surd$     &      & $\surd$     &  \\
	\cite{WangC.LiuF.WangJ.WeiW.Mei2016}& $\surd$     &      & $\surd$     &  \\
	\cite{Zhao2015}                     & $\surd$     &      & $\surd$     &  \\
	\cite{Zhao2016}                     & $\surd$     &      & $\surd$     &  \\
	\cite{Wei2015a}                     & $\surd$     &      & $\surd$     &  \\
	\cite{Wei2015e}                     & $\surd$     &      & $\surd$     &  \\
	\cite{Lannoye2012e}                 &      & $\surd$     & $\surd$     &  \\
	\cite{Lannoye2012d}                 &      & $\surd$     & $\surd$     &  \\
	\cite{Lannoye2012}                  &      & $\surd$     & $\surd$     &  \\
	\cite{Lannoye2015}                  &      & $\surd$     & $\surd$     &  \\
	\cite{Lu2005}                       &      & $\surd$     &      & $\surd$ \\
	\cite{Bresesti2003}                 &      & $\surd$     &      & $\surd$ \\
	\cite{Capasso2003}                  &      & $\surd$     &      & $\surd$ \\
    \hline
    \end{tabular}%
  \label{tab_metric}%
\end{table}%

\subsection{Physical Explanation of Grid-side Flexibility}
Operational flexibility has been defined as the ability of a power system to deploy its adjustable resources to respond to changes in net load, where net load is defined as the remaining system load not served by variable generation \cite{Lannoye2012e}. Note that there is no explicit definition of grid-side flexibility in the literature so far. Inspired by \cite{Lannoye2012e}, we define grid-side flexibility as below.

\begin{definition}
Grid-side flexibility is the ability of a power network to deploy its flexibility resources to cope with volatile changes of power system state in operation.
\end{definition}

The grid-side flexibility resources can be used to adjust the physical characteristics of power networks, such as topology or parameter, so as to enhance power system security and economy. According to the types of controllable variables, we classify grid-side flexibility resources into two categories: the discrete and the continuous. Discrete resources are capable of changing the topology of a power network, via transmission expansion planning (TEP)  and line switching (LS) in operation. Continuous resources, based on electronics controllers, which can offer control capability of power flow, voltage, phase angle and so on, include FACTS and HVDC.

\subsubsection{Discrete Flexibility Resources}
Typical discrete operational resources in grid side consist of line expansion planning and line switching, where the binary controllable variables indicate the availability of changing the ON/OFF state of lines.

Transmission expansion planning determines the time, the location, and the type of transmission lines to be built with the minimum investment cost, in order to guarantee the  adequacy of  supplying the forecasted electricity demand over the planning horizon. Traditionally, the electricity demand has been considered as a major source of uncertainty in planning decision-making. However, the increasing penetration of volatile renewable  generation contributes to additional uncertainty with a sustained growth, which has been becoming another main source of uncertainty \cite{Jabr2013a}. Different approaches have been adopted to reinforce grid-side flexibility in order to facilitate the integration of renewable  generation in a cost-effective way, such as robust optimization \cite{Jabr2013a,ChenB.WangJ.WangL.HeY.andWang2014}, stochastic optimization \cite{Park2013a,Yu2009a,Munoz2012a}.

Compared with the existing operational control methods such as generation re-dispatch or load shedding, line switching can control the topology of power networks, resulting in the change of power flow distribution, thus alleviating congestions. Line switching has been used for corrective action to mitigate transmission flow violations \cite{Bacher1986}, to improve power system security in contingency cases \cite{Schnyder1988,Schnyder1990}, and to manage transmission congestions \cite{Granelli2006}. The benefits of line switching to improve economic efficiency is first analyzed in \cite{O&apos;Neill2005} in a market environment. The impact of line switching on reducing production cost is investigated in \cite{Fisher2008}, and is extended in \cite{Hedman2008} to investigate the impacts of network topology changes on nodal prices, load payments, generation revenues, congestion costs, and flowgate prices. Line switching has also been introduced in security-constrained unit commitment in \cite{Khodaei2010b} for congestion management. In addition, line switching has been incorporated into chance-constrained \cite{Qiu2014} or adaptive robust model \cite{Taheri2015} to consider uncertainties. Whereas line switching enables controlling network topology, it is a discrete control method with limited control accuracy. On the other hand, N-1 contingency criterion is required to be satisfied when applying line switching, which imposing additional complexity in solving the related optimization problem.

\subsubsection{Continuous Flexibility Resources}
Typical continuous operational flexibility resources in grid side include FACTS and HVDC. Compared with generation-side and demand-side flexibility, the continuous resources of grid-side flexibility resources allow faster and cheaper controls to accommodate uncertainty. 

FACTS can be regulated very quickly and frequently, since the power electronics technology allows very fast response time (in the order of millisecond) \cite{Zhang2012a}. Moreover, FACTS are capable of controlling the interrelated parameters that govern the operation of transmission systems, including series impedance, shunt impedance, current, voltage, phase angle \cite{Hingorani2000a}. These advantages make it effective in handling the uncertainty of renewable  generation. Thyristor controlled series capacitor (TCSC) has been applied to make a better utilization of wind power generation in \cite{Yang2012,Nasri2014}. Phase-shifting transformer has shown excellent capability in mitigating congestion due to the forecast error of renewable  generation \cite{Thakurta2015}.

HVDC lines incur less power losses over long distances, while enhancing the controllability of the power grid, making it an attractive option to integrate renewable  generation \cite{Vrakopoulou2013}. There are mainly two types of HVDC techniques: line-commutated HVDC (LCC-HVDC) and voltage-source converter based HVDC (VSC-HVDC). The LCC-HVDC is the most popular and widespread HVDC technique nowadays \cite{Torres-Olguin2012}, while the emerging VSC-HVDC allows faster and reversible control of active and reactive power in a decoupling manner, facilitating the interconnection of AC power grid with offshore wind energy \cite{Fu2016}. Besides, VSC-HVDC can also be used to supply industrial networks and operate in weak transmission networks \cite{Cole2009}. Some significant contributions have been made to exploit the flexibility of HVDC with the integration of renewable energy generation. With regard to LCC-HVDC, a stochastic multi-period optimal power flow model is presented in \cite{Rabiee2014}, which includes an offshore wind farm connected to the AC power grid by a LCC-HVDC link. It is observed that the availability of transmission network capacity at the interface of AC/DC networks is a key factor affecting the utilization of wind power generation. As to VSC-HVDC, a security-constrained unit commitment solution is presented for the optimal integration of large-scale offshore wind energy into a power grid \cite{Fu2016}. The proposed multi-terminal HVDC grid interface provides a more flexible wind energy injection to the AC power grid, and mitigates inland AC network congestion. Probabilistic security-constrained optimal power flow algorithms incorporating renewable resources are proposed to take advantage of the HVDC line and meshed HVDC grids controllability for post-disturbance control actions \cite{Vrakopoulou2013,Vrakopoulou2013a}. 

\subsubsection{Combination of Discrete and Continuous Flexibility Resources}
A combination of discrete and continuous resources may be a promising way to fully exploit grid-side flexibility. An amount of effort has been devoted to incorporating continuous flexibility resources into transmission expansion planning. Phase-shifting transformer is considered as an element in transmission expansion planning to extend the utilization of classical components \cite{Miasaki2012a}. An investment valuation approach is proposed in a real option analysis framework to assess the option value of FACTS in transmission expansion planning \cite{Blanco2011a,Blanco2009a}. In \cite{Konstantelos2015a}, phase-shifting transformer is incorporated into transmission expansion planning based on a stochastic framework to investigate the potential for non-conventional assets to accommodate renewable energy. It is demonstrated that FACTS may reduce or defer the need for building new transmission lines \cite{Blanco2011a,Konstantelos2015a}. In addition, HVDC has also been incorporated in transmission planning problem for improving power system economic dispatch efficiency with wind energy integration \cite{Urquidez2015}, and facilitating the exploitation of offshore wind resources \cite{Torbaghan2015}.

\subsection{Mathematical Definition of Grid-side Flexibility}

Mathematically, grid-side flexibility resources provide additional operational  flexibility by introducing extra controllable variables into network constraints to control power flow. Inspired by \cite{Wu1982,Hou2015}, the mathematical definition of grid-side flexibility region is proposed. Let $\bm{P}=[\bm{P}_G,\bm{P}_L,\bm{P}_F]^{\top}$ denote the vector of active power injections, including power generation $\bm{P}_G$, load $\bm{P}_L$ and active power injection of FACTS or HVDC $\bm{P}_F$ with respect to power injection model \cite{Noroozian1997,Song2002}. Similarly, let $\bm{Q}=[\bm{Q}_G,\bm{Q}_L,\bm{Q}_F]^{\top}$ denote the vector of reactive power injections, where the subscripts have the same meanings as in the active power injection vectors. $\bm{z}$ denotes a vector of binary variables indicating the availability of lines. $\bm{V}$ and $\bm{\theta}$ denote voltage magnitude and angle, respectively. $\bm{f}$ and $\bm{h}$ represent all steady-state constraints, including economic dispatch (ED), transmission expansion planning and line switching. The extended steady-state security region $\Omega_{ESSR}$ is defined as below.

\begin{definition}
	The extended steady-state security region $\Omega_{ESSR}$ is defined by
\begin{equation}
\begin{aligned}
\Omega_{ESSR} := & \{(\bm{P},\bm{Q},\bm{z}) | \bm{f}(\bm{P},\bm{Q},\bm{z},\bm{V},\bm{\theta}) = \bm{0}, \\
& \bm{h}(\bm{P},\bm{Q},\bm{z},\bm{V},\bm{\theta})\le \bm{0}\}. \nonumber
\end{aligned}
\end{equation}
\end{definition}
When the variables related to grid-side flexibility resources (i.e. $\bm{P}_F$, $\bm{Q}_F$, $\bm{z}$) are fixed, and only the economic dispatch constraints are considered, $\Omega_{ESSR}$ is equivalent to the conventional steady-state security region \cite{Wu1982}.

The active power injection $\bm{P}$ can be re-ordered as $\bm{P}=[P_{Ref},\bm{P}_C,\bm{P}_U]^{\top}$ where $P_{Ref}$ is the active power injection of the reference bus, and $\bm{P}_C$ is the active power injection of all flexibility resources, including generator, load and FACTS/HVDC, while $\bm{P}_U$ is the uncertain active power injection of renewable energy generation. As a result, all flexibility resources can be denoted by $\bm{P}_{CR}=[P_{Ref},\bm{P}_C]^{\top}$. Similarly, $\bm{Q}$ is re-ordered as $\bm{Q}=[Q_{Ref},\bm{Q}_C,\bm{Q}_U]^{\top}$ where the subscripts have the same meaning as in the active power injection vectors. All controllable reactive power resources are denoted by $\bm{Q}_{CR}=[Q_{Ref},\bm{Q}_{C}]^{\top}$. Therefore, the admissible region \cite{WangC.LiuF.WangJ.WeiW.Mei2016} considering transmission control resources is defined as below.

\begin{definition}
$\Omega_{AR}$ is an admissible region on the subspace spanned by uncertain variable vectors $\bm{P}_{U}$ and $\bm{Q}_{U}$, if for any $(\bm{P}_{U},\bm{Q}_{U})$ in $\Omega_{AR}$, there exists at least one $(\bm{P}_{CR},\bm{Q}_{CR},\bm{z})$, so that $(\bm{P}_{CR},\bm{Q}_{CR},\bm{z}, \bm{P}_{U},\bm{Q}_{U})$ is in $\Omega_{ESSR}$. That is, 
\begin{equation*}
\begin{aligned}
\Omega_{AR}:= & \{(\bm{P}_{U},\bm{Q}_{U}) | \forall(\bm{P}_{U},\bm{Q}_{U})\in\Omega_{AR}, \\
& \exists(\bm{P}_{CR},\bm{Q}_{CR},\bm{z}), s.t. (\bm{P}_{CR},\bm{Q}_{CR},\bm{z},\bm{P}_{U},\bm{Q}_{U})\in\Omega_{ESSR}\}.
\end{aligned}
\end{equation*}
\end{definition}

The defined admissible region is such a region that every uncertain power injection in it can be accommodated by all the flexibility resources without violating steady-state constraints. Fixing the variables corresponding to grid-side flexibility resources results in another admissible region, $\Omega_{AR_0}$. Hence, the grid-side flexibility region is defined as below.

\begin{definition}
Grid-side flexibility region, $\Omega_{GR}$, is the additional region of $\Omega_{AR}$ when compared with the corresponding $\Omega_{AR_0}$, where grid-side flexibility is omitted. That is,
\begin{equation*}
\Omega_{GR} := \Omega_{AR} - \Omega_{AR_0}.
\end{equation*}
\end{definition}

The grid-side flexibility region denotes the expanded region of uncertain power injections where the uncertainty can be accommodated by grid-side flexibility resources. Since the grid-side flexibility region is defined on the subspace spanned by uncertain power injections, it may be discrete or continuous depending on the characteristic of $(\bm{P}_U,\bm{Q}_U)$. Besides, it is associated with the mixed-integer nonlinear non-convex constraints $\bm{f}$ and $\bm{h}$.

\subsection{Geometric Interpretation}\label{sec_example}
\begin{figure}[htb]
\centering
\includegraphics[width=3in]{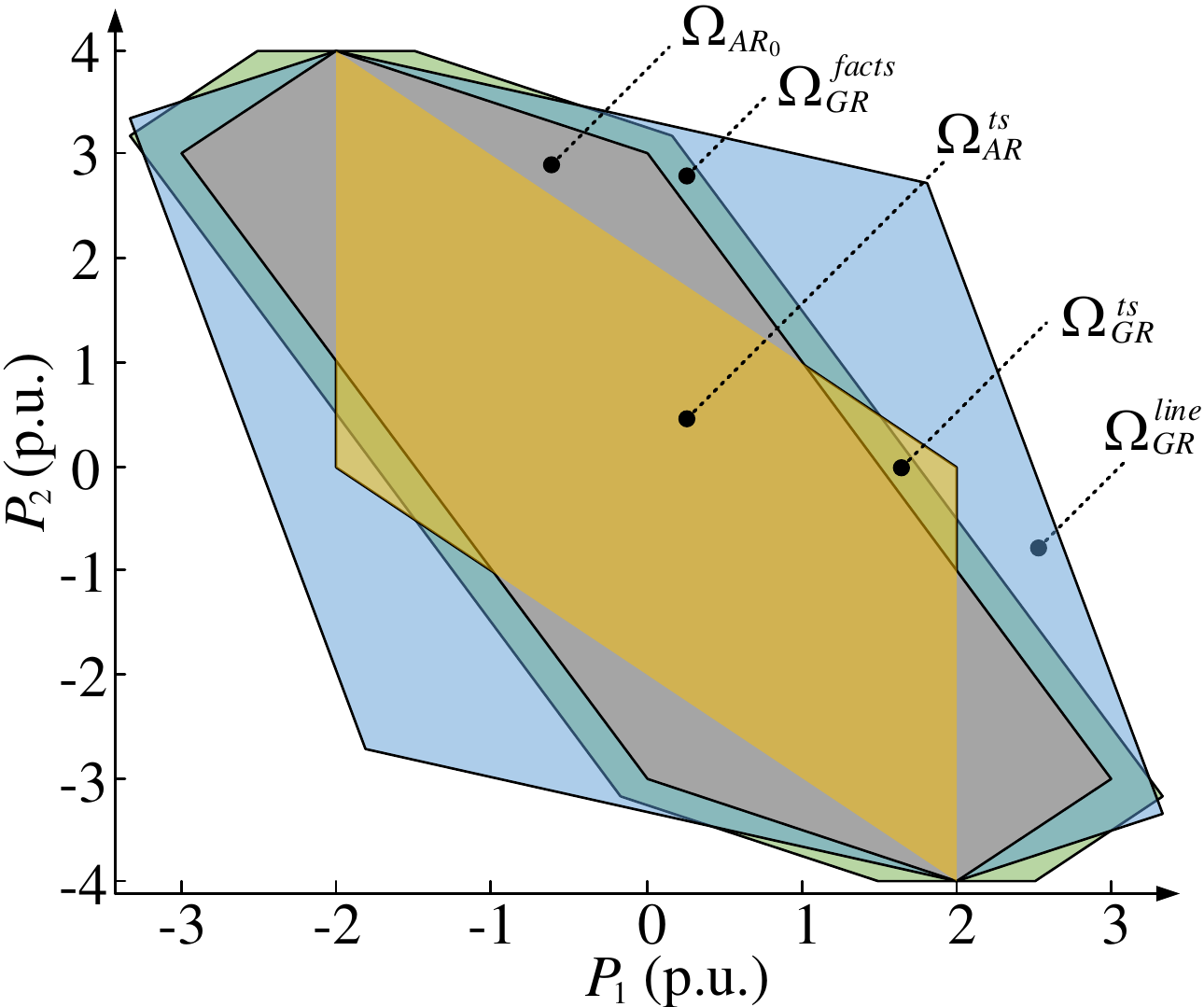}
\caption{Admissible regions with different grid-side flexibility resources.}
\label{fig_ar_trans}
\end{figure}

For simplicity and comprehensibility, an conceptual example of a modified 3-bus system  \cite{Hou2015} with DC power flow is employed. It is assumed that the admittance of each line is 1 p.u.. The lower and upper limits of line 1-2 and line 2-3 are -2 p.u. and 2 p.u., respectively, while the lower and upper limits of line 1-3 are -1 p.u. and 1 p.u., respectively. $P_1$ and $P_2$ are uncertain power injections of renewable energy generation at bus 1 and 2, respectively. $P_3$ is a controllable power injection at bus 3. The lower and upper limits of all the power injections are -4 p.u. and 4 p.u., respectively. The admissible regions can be constructed using Fourier-Motzkin elimination \cite{bertsimas1997introduction}. 
Figure. \ref{fig_ar_trans} depicts the admission regions with different grid-side flexibility
 resources. 

The extended steady-state security region is constructed as
\begin{equation}
\begin{aligned}
  \Omega_{ESSR} = & \left\{ 
  \begin{bmatrix}
  P_1 \\ P_2 \\ P_3
  \end{bmatrix} | 
  \begin{bmatrix}
  -2 \\ -2 \\ -1 \\ 0
  \end{bmatrix} \le \frac{1}{3} 
  \begin{bmatrix}
  1 & -1 & 0 \\
  1 &  2 & 0 \\
  2 &  1 & 0 \\
  1 &  1 & 1
  \end{bmatrix}
  \begin{bmatrix}
  P_1 \\ P_2 \\ P_3
  \end{bmatrix} \le
  \begin{bmatrix}
  2 \\ 2 \\ 1 \\ 0
  \end{bmatrix}, \right. \\
  & \left. \begin{bmatrix}
  -4 \\ -4 \\ -4
  \end{bmatrix} \le
  \begin{bmatrix}
  P_1 \\ P_2 \\ P_3
  \end{bmatrix} \le
  \begin{bmatrix}
  4 \\ 4 \\ 4
  \end{bmatrix} \right\}.
\end{aligned}
\end{equation}

Without transmission control, $\Omega_{AR_0}$ is derived as
\begin{equation}\label{eq_arp}
  \begin{bmatrix}
  -6 \\ -6 \\ -3
  \end{bmatrix} \le
  \begin{bmatrix}
  1 & -1 \\
  1 &  2 \\
  2 &  1
  \end{bmatrix}
  \begin{bmatrix}
  P_1 \\ P_2
  \end{bmatrix} \le
  \begin{bmatrix}
  6 \\ 6 \\ 3
  \end{bmatrix}
\end{equation} 
which represents the gray hexagonal region in Fig. \ref{fig_ar_trans}.

(1) \emph{Installation of a FACTS device.} Suppose a series FACTS device with the capacity of 0.5 p.u. is installed on line 1-3. According to the power injection model, the device is formulated as two active power injections at buses 1 and 3, respectively. 
Hence, $\Omega_{AR}$ is constructed as
\begin{equation}
  \begin{bmatrix}
  -4 \\ -9 \\-6.5 \\ -6.5 \\ -3.5
  \end{bmatrix} \le 
  \begin{bmatrix}
  0 & 1 \\
  1 & 2 \\
  1 & -1 \\
  1 & 2 \\
  2 & 1
  \end{bmatrix}
  \begin{bmatrix}
  P_1 \\ P_2
  \end{bmatrix} \le
  \begin{bmatrix}
  4 \\ 9 \\ 6.5 \\6.5 \\ 3.5
  \end{bmatrix}.
\end{equation}
Given any $P_2$, by regulating the power flow through line 1-3, FACTS device expands the admissible range of $P_1$, yielding the extra region $\Omega_{GR}^{facts}$.

(2) \emph{Construction of a new line.} Suppose the transmission capacity of corridor 1-3 is doubled by building a new line. Thus, the admittance of the corridor 1-3 is reduced to 0.5 p.u.. $\Omega_{GR}^{line}$ is constructed as below.
\begin{equation}
  \begin{bmatrix}
  -10 \\ -10 \\ -10
  \end{bmatrix} \le 
  \begin{bmatrix}
  1 & -2 \\
  1 &  3 \\
  4 &  1
  \end{bmatrix}
  \begin{bmatrix}
  P_1 \\ P_2
  \end{bmatrix} \le
  \begin{bmatrix}
  10 \\ 10 \\ 10
  \end{bmatrix}
\end{equation}
Compared with FACTS, the new transmission line expands the admissible range of $P_1$ further. However, the expansion is within a little smaller range of $P_2$.

(3) \emph{Line switching.}  Suppose line switching is applied in line 1-3. When line 1-3 is switched on, $\Omega_{AR_0}$ is the same as (\ref{eq_arp}), whereas when line 1-3 is switched off, $\Omega_{AR_0}^{ts}$ is identified as
\begin{equation}\label{eq_ts_off}
  \begin{bmatrix}
  -2 \\ -2
  \end{bmatrix} \le 
  \begin{bmatrix}
  1 & 0 \\
  1 & 1
  \end{bmatrix}
  \begin{bmatrix}
  P_1 \\ P_2
  \end{bmatrix} \le
  \begin{bmatrix}
  2 \\ 2
  \end{bmatrix}
\end{equation}
which represents the parallelogram region in Fig. \ref{fig_ar_trans}. Thus, the complete $\Omega_{AR_0}$ is the union of the regions defined by (\ref{eq_arp}) and (\ref{eq_ts_off}), yielding a concave region. It is observed that line switching affects the convexity of admissible region, and the resulting two triangle grid-side flexibility regions, $\Omega_{GR}^{ts}$, depend on system parameters and operational state. 

\section{Problem Formulations Corresponding to Grid-side Flexibility} \label{sec_form}
Depending on the formulation of uncertain variables, the models considering the grid-side flexibility resources can be classified into two categories: stochastic programming (SP) and robust optimization (RO). Stochastic programming is based on the assumption that the complete information about the probability distribution of uncertain variables is available and usually is focused on the expectation of total cost. Robust optimization does not explicitly require the probability distribution of uncertain variables, and typically uses uncertainty sets to characterize uncertainty. 

\subsection{Stochastic Programming (SP)}\label{sec_form_sp}
Taking advantage of the probability distribution of random variables, stochastic programming usually aims to minimize the expected total cost. Different stochastic programming approaches associated with grid-side flexibility can be categorized as follows.

\subsubsection{Scenario-Based Model}
Scenario-based approach uses discrete scenarios to simulate the possible realization of uncertainty. With a presumed probability distribution, a large number of scenarios are generated to achieve an acceptable solution accuracy. However, more scenarios lead to a heavier computational burden, and even intractability in a large-scale system. Thus an important issue is to weigh the trade-off between accuracy and tractability. In \cite{Nasri2014}, wind power generation uncertainty is modeled using a set of plausible scenarios and the model is cast as a two-stage stochastic programming problem to minimize wind power spillage with FACTS. TCSC and fixed series capacitor (FSC) are incorporated into the AC optimal power flow (OPF) model to optimally determine the setting of such devices. A stochastic multi-period OPF model is presented in \cite{Rabiee2014} which contains an offshore wind farm connected to the grid by a LCC-HVDC link. The uncertainty of wind power generation is characterized by a scenario-based approach \cite{Atwa2011}. With respect to transmission expansion planning, different scenarios of load levels and wind power generation are employed to capture the uncertainty \cite{Munoz2012a,Orfanos2013a,Villumsen2013a}, and a scenario tree is used to describe the evolution of wind power generation capacity \cite{Konstantelos2015a}. Additionally, approximate partitions are employed to deal with dependent random variables including wind availability and load \cite{Park2013a}.

\subsubsection{Point-Estimation-Based Model}
Point estimation method provides an approximate description of the statistical information of the functions associated with random variables. Two-point estimation method has been applied in transmission expansion planning with the presence of wind power generation \cite{Moeini-Aghtaie2012a} and under smart grid environment with demand response resources \cite{Hajebrahimi2015a}. To avoid the computational burden of Monte Carlo simulation and mathematical calculations associated with analytical methods, approximate point estimation (2m+1)-scheme is used in \cite{Arabali2014a} to perform probabilistic OPF for transmission planning.

\subsubsection{Chance-Constrained Programming Model}
Chance-constrained program allows the constraints to be violated with a predefined small probability. A novel approach is proposed in \cite{GuhaThakurta2015} to quantify the amount of achievable confidence of system operation and its improvement with the help of power flow controlling devices (PFCs), such as phase-shifting transformers and HVDC. The influence of PFCs is considered using the generalized power transfer distribution factor, and the chance constraints on line flows are transformed into nonlinear deterministic constraints. A security-constrained OPF algorithm is proposed in \cite{Vrakopoulou2013a}, which incorporates the stochastic infeed of renewable resources using chance constraints, and takes advantage of the HVDC line controllability. Each HVDC line is approximated by two virtual voltage sources located at the two nodes where the HVDC line is connected. The formulation in \cite{Vrakopoulou2013a} is extended in \cite{Vrakopoulou2013} to investigate multi-terminal HVDC (MTDC) grids. A chance-constrained model is proposed in \cite{Qiu2014} to explore the possibility of changing transmission topology, by line switching, to accommodate higher penetration of wind power generation and reduce the fuel costs of thermal units. Chance-constrained programming formulation has also been employed to tackle the uncertainties of load and wind power generation in transmission expansion planning \cite{Yu2009a}.

\subsection{Robust Optimization (RO)}\label{sec_form_ro}
Instead of the requirement of probability distribution information in stochastic programming, robust optimization requires moderate information about uncertainty, such as the mean and the range of the uncertain variables. In addition, the optimal solution immunizes against all possible realizations of uncertain variables within a deterministic uncertainty set. However, this approach may be over-conservative, since the worst case rarely occurs. Different robust optimization frameworks associated with grid-side flexibility are categorized as follows.

\subsubsection{Uncertainty-Set-Based Model}
Characterizing the uncertain renewable energy generation by polyhedral uncertainty sets, these single-stage \cite{Jadid2013a} or two-stage \cite{Jabr2013a,ChenB.WangJ.WangL.HeY.andWang2014,Ruiz2015} frameworks minimize the costs or regret \cite{ChenB.WangJ.WangL.HeY.andWang2014} in the worst-case scenario. The frameworks have been applied not only in transmission expansion planning \cite{Jadid2013a,Jabr2013a,ChenB.WangJ.WangL.HeY.andWang2014,Ruiz2015}, but also in line switching \cite{Taheri2015}. In \cite{Taheri2015}, an adaptive robust optimal line switching approach is proposed to minimize the total operational and penalty costs pertaining to the worst-case scenario of uncertain net nodal load.

\subsubsection{Taguchi's Orthogonal-Array-Testing-Based Model}
Taguchi's orthogonal array testing (TOAT) selects representative scenarios to approximate the uncertainty space such that only a small number of tests are  performed without sacrificing too much accuracy. The selected worst-case uncertainty scenarios are the orthogonal corners of the approximate uncertainty space \cite{Kazemi2014a}. A robust transmission expansion planning method is proposed based on TOAT in \cite{Yu2011}, which considers the uncertainty of load and renewables.

\subsubsection{Information-Gap-Decision-Theory-Based Model}
This kind of method uses information-gap decision theory (IGDT) to find the robust feasibility region against all uncertain variables by bounding the allowable range of objective function regarding a definite uncertainty budget \cite{Kazemi2014a}. In \cite{Kazemi2014a}, an IGDT-based framework is presented to cope with the uncertain nature of capital costs and electricity demands within transmission expansion planning. The uncertainty is formulated by a distribution-free envelope-bound model. Compared with other robust optimization frameworks, this IGDT model can explicitly represent the robustness cost.

\section{Solution Approaches }
In this section, we review the solution approaches which are applied to solve the problems in Section \ref{sec_form}.

\subsection{Stochastic Programing}
The solution approaches pertaining to the three categories of formulations in Section \ref{sec_form_sp} are summarized as follows.

\subsubsection{Scenario-Based Model}
When formulated as a nonlinear programming (NLP) problem \cite{Nasri2014}, \cite{Rabiee2014} or a mixed-integer linear programming (MILP) problem \cite{Munoz2012a}, the stochastic programming problem may be solved by  commercial solvers, such as CONOPT, CPLEX. Decomposition techniques are usually employed to improve the computational efficiency. In \cite{Orfanos2013a}, a Benders decomposition algorithm is adopted to tackle the proposed stochastic transmission expansion planning problem. A Benders cut is generated and added cumulatively into the master problem when any constraint is active in any of the operation subproblems. Benders decomposition is also implemented in \cite{Konstantelos2015a} to decompose the original problem into a master problem of transmission investment and several subproblems of power system operation. Unlike Benders decomposition to add additional constraints in the master problem, the Dantzig-Wolfe decomposition adds extra columns (variables) in the master problem. This method is used to solve the line capacity expansion problem with line switching under future uncertainty in demand and wind generation capacity \cite{Villumsen2013a}. In addition, a sequential approximation approach is applied to solve the two-stage stochastic transmission expansion planning model in \cite{Park2013a}, which is an iterative method that solves a two-stage stochastic problem at each iteration.

\subsubsection{Point-Estimation-Method-Based Model}
The non-dominated sorting genetic algorithm (NSGA II) is a type of evolutionary algorithms, which is capable of handling the nonlinearity and mixed-integer nature of stochastic transmission expansion planning with point estimation method \cite{Moeini-Aghtaie2012a,Hajebrahimi2015a,Arabali2014a}. In NSGA II, the first population is initialized and sorted into different Pareto fronts. The Pareto fronts and their individuals are ranked according to their non-dominancy levels. A binary tournament algorithm based on the non-dominancy rank and crowding distance is used to select the parent populations. Then, traditional crossover and mutation operators are adopted to generate offspring populations. Finally, the parents and offspring are composed to form a collection where a next generation is selected. This process continues until the termination criterion is satisfied.

\subsubsection{Chance-Constrained Model}
Chance-constrained problems can be solved using different approaches. In \cite{GuhaThakurta2015}, the chance constraint is transformed into a nonlinear deterministic constraint under the assumption that the random variable follows a normal distribution. The reformulated problem is solved using the KNITRO solver. An approximation approach named sample average approximation (SAA) is applied to solve the chance-constrained line switching problem in \cite{Qiu2014}. Sample scenarios are drawn from the distribution. The frequency of the scenarios, where the constraints are satisfied, is used to approximate the probability that the constraints hold. Another scenario approach is proposed in \cite{Vrakopoulou2013a,Vrakopoulou2013} to address the chance-constrained problems. The scenario-based approach is first used to determine the minimum volume set with a desired probability level, then a robust problem is formulated where the uncertainty is confined in the set. Additionally, a two-step genetic algorithm is designed to solve the chance-constrained transmission expansion problem in \cite{Yu2009a}. In the first step, a few genetic algorithm optimization generations are executed and a normal distribution is adopted to simulate wind power output. In the second step, the real wind power output distribution model is adopted.

\subsection{Robust Optimization}
Corresponding to the formulations in Section \ref{sec_form_ro}, the solution approaches are reviewed in three categories:

\subsubsection{Uncertainty-Set-Based Model}
A single-stage robust optimization problem can be formulated as an MILP problem \cite{Jadid2013a} and solved by commercial solvers, such as CPLEX. The widely used two-stage robust problem can be solved using decomposition techniques or meta-heuristic algorithms. When applying decomposition techniques, the robust optimization is usually formulated as a master problem and a subproblem, where the Big-M method is usually used to linearize the nonlinear terms through a set of auxiliary binary/continuous variables and disjunctive constraints \cite{Taheri2015}. One of the typical decomposition techniques is the Benders decomposition, which is based on the dual information of the second-stage problem. A Benders decomposition scheme is studied in \cite{Jabr2013a}, where the algorithm iterates between a master problem that minimizes the cost of the transmission expansion plan and a subproblem that minimizes the maximum curtailment of load and renewable generation. Different from Benders decomposition, the column-and-constraint generation (C\&CG) \cite{Zeng2013} technique involves primal cuts rather than dual cuts, which performs computationally better than a Benders decomposition algorithm \cite{Ruiz2015}. This technique has been used to solve robust transmission expansion planning problem \cite{ChenB.WangJ.WangL.HeY.andWang2014,Ruiz2015} and optimal line switching problem \cite{Taheri2015}. Note that the primal cuts add new variables and constraints into the master problem at every iteration, while dual cuts only add new constraints, which means the primal cuts increase the size of the master problem more than dual cuts do at every iteration. This feature may lead to intractability of the master problem after too many iterations. Meta-heuristic method is another technique to solve a robust optimization problem. In \cite{Chatthaworn2015}, the adaptive Tabu search is utilized to solve the proposed robust transmission expansion planning problem, where AC power flow constraints are enforced. In \cite{Wen2015}, a hybrid algorithm combing greedy search and particle swarm optimization (PSO) is adopted, where greedy search is used for local search and PSO is used for global search.

\subsubsection{Taguchi's Orthogonal-Array-Testing-Based Model}
A genetic algorithm is applied to solve a robust transmission expansion planning problem based on Taguchi's orthogonal array testing \cite{Yu2011}. In order to incorporate the information of all the scenarios in the initial genetic algorithm population, the population is divided into some subsets which are initialized by different testing scenarios. Besides, a scheme set which stores the already checked planning schemes and the checked results is proposed to avoid repeatedly checking the same candidate scheme in different random generations during the optimization process.

\subsubsection{Information-Gap-Decision-Theory-Based Model}
The augmented $\varepsilon$-constraint method is employed to solve the robust transmission expansion planning problem using information-gap decision theory in \cite{Kazemi2014a}. In order to overcome the drawbacks of the conventional $\varepsilon$-constraint method, lexicographic optimization is adopted to calculate the effective ranges for all objective functions. Moreover, augmented $\varepsilon$-constraint method is used to ensure the efficiency of Pareto optimal solutions, which optimizes a main augmented objective function chosen among all objective functions.

Table \ref{tab_form} summarizes the aforementioned formulations and solution approaches related to grid-side flexibility.

\begin{table*}[htbp]
  \centering
  \tiny
  \caption{Formulations and Solution Approaches of the Studies Associated with Grid-side Flexibility}
    \begin{tabular}{cccccccccc}
    \hline
    Ref. No. & \multicolumn{4}{c}{Grid-side flexibility Resource} & \multicolumn{2}{c}{Model Formulation} & \multicolumn{2}{c}{Power Flow Formulation} & Solution Approach \\
          & FACTS & HVDC  & TS    & TEP   & SP    & RO    & AC    & DC    &  \\
    \hline
	\cite{Nasri2014}                        & $\surd$     &      &      &      & $\surd$     &      & $\surd$     &      & NLP \\   	 
	\cite{GuhaThakurta2015}                 & $\surd$     & $\surd$     &      &      & $\surd$     &      &      & $\surd$     & NLP \\
	\cite{Rabiee2014}                       &      & $\surd$     &      &      & $\surd$     &      & $\surd$     &      & NLP \\
	\cite{Vrakopoulou2013a}                 &      & $\surd$     &      &      & $\surd$     &      &      & $\surd$     & Scenario Approach \\
	\cite{Vrakopoulou2013}                  &      & $\surd$     &      &      & $\surd$     &      &      & $\surd$     & Scenario Approach \\
	\cite{Qiu2014}                          &      &      & $\surd$     &      & $\surd$     &      &      & $\surd$     & Sample Average Approximation \\
	\cite{Yu2009a}                          &      &      &     & $\surd$      & $\surd$     &      &      & $\surd$     & Genetic Algorithm \\
	\cite{Park2013a}                        &      &      &     & $\surd$      & $\surd$     &      &      & $\surd$     & Sequential Approximation Approach \\
	\cite{Munoz2012a}                       &      &      &     & $\surd$      & $\surd$     &      &      & $\surd$     & MILP \\
	\cite{Moeini-Aghtaie2012a}              &      &      &     & $\surd$      & $\surd$     &      &      & $\surd$     & NSGA-II \\
	\cite{Orfanos2013a}                     &      &      &     & $\surd$      & $\surd$     &      &      & $\surd$     & Benders Decomposition\\
	\cite{Hajebrahimi2015a}                 &      &      &     & $\surd$      & $\surd$     &      &      & $\surd$     & NSGA-II \\
	\cite{Arabali2014a}                     &      &      &     & $\surd$      & $\surd$     &      &      & $\surd$     & NSGA-II \\
	\cite{Konstantelos2015a}                & $\surd$      &      &     & $\surd$      & $\surd$     &      &      & $\surd$     & Benders Decomposition \\
	\cite{Villumsen2013a}                   &      &      & $\surd$     & $\surd$      & $\surd$     &      &      & $\surd$     & Dantzig-Wolfe Decomposition \\
	\cite{Taheri2015}                       &      &      & $\surd$     &      &      & $\surd$     &      & $\surd$     & C\&CG \\
	\cite{Jabr2013a}                        &      &      &      & $\surd$     &      & $\surd$     &      & $\surd$     & Benders Decomposition \\
	\cite{ChenB.WangJ.WangL.HeY.andWang2014}&      &      &      & $\surd$     &      & $\surd$     &      & $\surd$     & C\&CG \\
	\cite{Ruiz2015}                         &      &      &      & $\surd$     &      & $\surd$     &      & $\surd$     & C\&CG \\
	\cite{Jadid2013a}                       &      &      &      & $\surd$     &      & $\surd$     &      & $\surd$     & MILP \\
	\cite{Kazemi2014a}                      &      &      &      & $\surd$     &      & $\surd$     &      & $\surd$     & Augmented $\varepsilon$-constraint method \\
	\cite{Yu2011}                           &      &      &      & $\surd$     &      & $\surd$     &      & $\surd$     & Genetic Algorithm \\
	\cite{Chatthaworn2015}                  &      &      &      & $\surd$     &      & $\surd$     & $\surd$     &      & Adaptive Tabu Search \\
	\cite{Wen2015}                          &      &      &      & $\surd$     &      & $\surd$     &      & $\surd$     & Greedy Search, PSO \\
    \hline
    \end{tabular}%
  \label{tab_form}%
\end{table*}%

\section{Future Research Directions and Challenges}

In this section, we discuss some future research directions and challenges of exploiting grid-side flexibility in power system operation and planning.

\subsection{Consideration of Power Network Topology and Parameter Changes}

Due to line switching, change of transformer tap, or N-k security criterion, the discrete changes of power network topologies and parameters during operation may impact power system economics and security. In order to formulate these discrete changes, integer variables are inevitably introduced into the model, making it an NP-hard problem. One approach to solve this problem is to cast it as an MILP problem, so that some commercial solvers can be adopted. Some approximations have to be made for linearizion, which may ignore the reactive power constraints. Moreover, there is still lack of high-efficiency algorithm to solve a two-stage robust optimization problem with integer variables introduced in both the first and second stages,
for instance, a transmission expansion planning problem with unit commitment in the second stage, or a transmission expansion planning problem with line switching in the second stage, since the global optimal solution of the second stage may not be guaranteed in general. Most studies in the literature only consider integer variables in the first stage \cite{Khodaei2010a,Dehghan2015,Chen2016}, limiting the underlying flexibility in the second stage. When the scale of the second stage problem is small, enumeration algorithm may be employed.
Dantzig-Wolfe decomposition and a branch-and-price algorithm have been developed to solve a two-stage stochastic optimization problem with integer variables in both stages \cite{Villumsen2013a}, which may be extended to solve a robust optimization problem.
A nested column-and-constraint generation algorithm is presented in \cite{Zhao2011} to solve a linear two-stage robust optimization model with a mixed integer recourse problem.

\subsection{Consideration of Impacts of Reactive Power Flow}

Grid-side flexibility resources can not only change active power flow, but also control reactive power flow and improve voltage profile by generating/absorbing  extra reactive power. DC power flow formulation has been widely used to incorporate grid-side flexibility in power system operation and planning, which requires much less computational burdens and is appropriate in the case that reactive power flow can be neglected. However, in order to make full use of grid-side flexibility, the impacts on the reactive power flow need to be taken into account. In other words, AC power flow constraints are required to be considered in the optimization problem, resulting in an AC optimal power flow (AC OPF) problem, which is difficult to address due to the inherent nonlinearity and non-convexity. Moreover, the wide controllable range of FACTS parameter makes the problem even harder to converge.

In the literature, some references have dedicated to solve this difficult problem. Stochastic AC OPF models with FACTS devices \cite{Nasri2014} or HVDC links \cite{Rabiee2014} have been proposed, which are solved using CONOPT solver. The impacts and benefits of unified power flow controller (UPFC) to wind power integration in unit commitment are studied using a two-stage stochastic model with AC constraints \cite{Li2018}. A robust transmission planning model with AC power flow constraints has been presented in \cite{Chatthaworn2015} and solved by an adaptive Tabu search algorithm. Under the deterministic assumption, some other approaches have been developed, such as interior-point method \cite{Song2002,Zhang2001a,Zhang2001b,Zhang2012}, sequential quadratic programming \cite{Lehmkoster2002}, conic programming \cite{Garcia-bellido2012}, linear approximation \cite{Akbari2014a}, meta-heuristic algorithm \cite{Panuganti2013,Chung2001,Bhasaputra2002,Leung2011}.
Particularly, in radial AC networks, semi-definite programming (SDP)\cite{Lavaei2012} and second-order cone programming (SOCP) \cite{Jabr2006,Sojoudi2012,Farivar2013a,Low2014,Low2013,Huang2016} are adopted for convexification of OPF problem.
Despite all these significant contributions, it is still necessary to develop efficient algorithms to fully exploit grid-side flexibility in AC networks regardless of topology.

\subsection{Market for Grid-side Flexibility}
It is argued that the current market compensation mechanism does not provide incentives for the operation of FACTS devices in a way that maximizes the system benefit, since the owner of the devices only receives a regulated rate of return after installation, and it is against the owner’s interest to change the set-point more frequently to avoid additional stress on the device and higher maintenance \cite{Sahraei-Ardakani2015}. Therefore, an efficient market design and better compensation mechanism for grid-side flexibility is necessary. A real-time electricity market design is presented in \cite{ONeill2008}, where all assets including generation, load, and transmission are allowed to bid and transmission lines are compensated for both capacity and admittance, providing incentives for efficient operation of transmission-related assets. However, it is pointed out that such a market design induces a positive externality problem, where an FACTS device may get rewarded because of the actions taken by another device \cite{Sahraei-Ardakani2015,Sahraei-Ardakani2012}. A sensitivity-based method is developed in \cite{Sahraei-Ardakani2015} to avoid the positive externality characteristics and provide economic incentives for the device owners to operate them in a socially optimal way. It should be noted that reactive power flow is not taken into account and the possibility exists that the sensitivity calculation may be unreliable.

Financial transmission rights (FTRs) were initially created to allow transmission users to hedge congestion cost \cite{ONeill2008}. Two types of such rights have been proposed: flowgate rights and point-to-point rights. Flowgate rights entitle holders to the dispatch shadow price for specific transmission assets or groups of assets, known as flowgates \cite{Chao2000,Chao1996}. Point-to-point rights entitle holders to settlement equal to the locational marginal price difference between the two points times the megawatt amount of the right \cite{Hogan1992}. Series FACTS devices have been embedded in the FTR auction through the power injection model to make more existing transmission capabilities to be allocated to market participants as FTRs \cite{WangX.2002}. A drawback in using FTR-style rights for compensating merchant transmission projects is that it is difficult to identify which particular set of FTRs would correspond to a transmission expansion project, and the order in which projects are built affects the rights awarded \cite{ONeill2008}.
There are some interesting topics remained to be investigated:

(1) A versatile compensation scheme for different grid-side flexibility resources. Previous studies mainly focused on the continuous resources, such as FACTS. Note that the discrete resources, such as line switching, may also benefit congestion management, highlighting the need of a generalized compensation mechanism for different controllable resources.

(2) FTR bidding strategy considering the uncertainty of renewable energy generation. The inherent uncertainty of renewable energy generation makes it more difficult to forecast congestions in power system. Therefore, an appropriate bidding strategy under uncertainty may make more profits and avoid possible high risk.

(3) Measuring and mitigating the market power of grid-side flexibility resources. The metrics of market power need to be investigated, since grid-side flexibility resources may change the topology and parameter of a power network to obtain market power. The problem may be resolved similarly to generator market power, through enforcement of market power mitigation measures \cite{ONeill2008}.

\subsection{Coordination of Different Flexibility Resources}
Enhancing the deliverability of generation-side and demand-side flexibility, grid-side flexibility is an efficient supplement to generation-side and demand-side flexibility. Coordinating the different resources of flexibility may help to improve power system security and economics under uncertainty, facilitating the integration of renewable energy generation with less investment. Some contributions have been made in the coordination of generation-side and grid-side flexibility in economic dispatch \cite{Yang2012,Nasri2014,Thakurta2015,Song2002,GuhaThakurta2015} system expansion planning \cite{Hemmati2013d,Roh2009a,Tor2008a,Khodaei2012a,Pozo2013a,Heidari2015a,AlvarezLopez2007a}, transmission-generation-demand response co-optimization \cite{Ozdemir2015a}, transmission expansion planning with demand response resources \cite{Hajebrahimi2015a}.

However, there are two main issues that have not been fully addressed yet. One is the coordination of multiple time scales of different flexibility resources. For instance, transmission expansion planning usually corresponds to the time scale of years, while the flexibility resources during operation are related to the time scale of hours, rendering a multi-stage time-coupled optimization problem, which is difficult to solve. The other is the potential interactions among different sources of flexibility. Therefore, the constraints on the interactions need to be added into the model, increasing the scale and complexity. A single intertemporal optimization framework is proposed in \cite{Haller2012} where long-term transition processes are driven by $\text{CO}_2$ prices, endogenous technological learning and increasing fuel costs, while characteristic time slices are used to represent short-term temporal fluctuations of supply and demand. However, the conceptual model can only provide qualitative results.

\subsection{Co-optimization of Multiple Energy Systems}
Nowadays, as is often the case, different energy systems are operated independently. However, synergy effects among various energy carriers can be achieved by taking advantage of their specific virtues \cite{GeidlM.KoeppelG.Favre-PerrodP.KlocklB.AnderssonG.FrohlichK.2007}. Electricity, for instance, can be transmitted quickly over long distances; natural gas can be stored using relatively simple and cheap technologies. Combing multiple energy systems may improve system efficiency and enhance reliability. 

The concept of energy hub was proposed in \cite{GeidlM.KoeppelG.Favre-PerrodP.KlocklB.AnderssonG.FrohlichK.2007}, which can be identified as a unit that provides the basic features in- and output, conversion, and storage of different energy carriers, such as electricity, natural gas and heat \cite{Geidl2007}. Effort has been made on utilizing energy hub, including optimal power flow \cite{Geidl2007}, interconnector \cite{Kienzle2008} and planning \cite{Kienzle2009,Zhang2015b,Zhang2015a}. Natural gas network has been introduced in power system operation \cite{Li2008,Liu2009b,Sahin2011,Liu2010} and planning \cite{Zhang2015b} in recent years. The effectiveness of firming variable wind energy generation by coordinating electricity and natural gas systems has been shown in \cite{Alabdulwahab2015}. Utilization of combined heat and power (CHP) unit \cite{Chen2015a} and the heat storage capacity of heating network \cite{Li2016,Li2016b} has been proven to be a cost-effective method to enhance power system flexibility to accommodate renewable energy generation. In addition, exploiting the flexibility of other forms of energy, such as hydro \cite{Khodayar2013,Jiang2012a,Abreu2012}, electric vehicles \cite{Shao2014,Ghofrani2014,Shao2015}, has also been reported to help to promote the integration of renewable energy.

As a result, the concept of grid-side flexibility may be extended to a multi-energy system, incorporating the network constraints of different energy systems. However, there are some challenging issues that have not been completely addressed. First, the combination of multiple energy systems introduces a massive quantity of nonlinear constraints, yielding a non-convex problem which is hard to solve. Certain linear approximations may be applied to simplify the problem. However, the trade-off between accuracy and tractability should be made carefully, and the applicability of the approximation should be specified. Second, different energy systems have distinct operational time scales due to inherent physical characteristics. For example, electricity is transmitted much faster than natural gas or heat. Therefore, the coordination of different time scales and dynamic processes may need to be considered, resulting in a complicated model with both algebraic equations and differential equations.

\section{Conclusions}
In this paper, we review the metrics corresponding to power system flexibility, and introduce the novel concept of grid-side flexibility in both physical and mathematical ways. Moreover, conceptual examples clearly demonstrate the impacts of grid-side flexibility in a geometric manner. Two major types of problem formulations with respect to grid-side flexibility are reviewed, and the associated solution approaches are surveyed. Some challenging problems nowadays and possible future research directions are discussed.

Grid-side flexibility is an effective supplement to generation-side and demand-side flexibility to accommodate the uncertainty of renewable energy generation. In order to fully utilize grid-side flexibility, topology changes and reactive power flow need to be considered, resulting in mixed-integer nonlinear optimization problems. Thus, efficient algorithms are required to solve the complicated problems. Efficient market mechanism is also required to provide economic incentives for optimal operation and investment of grid-side flexibility resources. Additionally, the co-optimization of different resources of flexibility and multiple energy systems may lead to a better way of coping with large-scale volatile renewable generations.

\section*{Acknowledgments}
This work was supported by the Research on Green Dispatch Technology of Shanxi Power Grid and the Foundation for Innovative Research Groups of the National Natural Science Foundation of China (51621065).

\section*{References}

\bibliography{library}

\end{document}